\documentclass[preprint,superscriptaddress]{revtex4}

\usepackage{graphicx}
\usepackage{array}
\usepackage{amssymb}
\usepackage{amsfonts}
\usepackage{amsmath}
\usepackage{mathrsfs}
\usepackage{color}
\usepackage{booktabs}
\usepackage{threeparttable}
\usepackage{multirow}
\usepackage{subfigure}
\usepackage{times}
\usepackage{epsfig}
\usepackage{threeparttable}
\usepackage{chngpage}
\usepackage{float}
\usepackage{algorithm}
\usepackage{algorithmic}
\usepackage{color}

\makeatletter
\renewcommand\@biblabel[1]{#1.}
\makeatother

\begin{document}

\title{Reduce the rank calculation of a high-dimensional sparse matrix based on network controllability theory}

\author{Chen Zhao}
\affiliation{College of Computer and Cyber Security, Hebei Normal University, Shijiazhuang, 050024, P. R. China\\}
\affiliation{Hebei Provincial Engineering Research Center for Supply Chain Big Data Analytics \& Data Security, Shijiazhuang, 050024, P. R. China\\}
\affiliation{Hebei Key Laboratory of Network and Information Security, Shijiazhuang 050024, P. R. China\\}

\author{Yuqing Liu}
\affiliation{College of Computer and Cyber Security, Hebei Normal University, Shijiazhuang, 050024, P. R. China\\}
\affiliation{Hebei Provincial Engineering Research Center for Supply Chain Big Data Analytics \& Data Security, Shijiazhuang, 050024, P. R. China\\}

\author{Li Hu}
\affiliation{College of Computer and Cyber Security, Hebei Normal University, Shijiazhuang, 050024, P. R. China\\}
\affiliation{Hebei Key Laboratory of Network and Information Security, Shijiazhuang 050024, P. R. China\\}

\author{Zhengzhong Yuan}\email{zyuan@mnnu.edu.cn}
\affiliation{School of Mathematics and Statistics, Minnan Normal University, Zhangzhou, 363000, P. R. China\\}
\affiliation{Fujian Key Laboratory of Data Science and Statistics, Minnan Normal University, Zhangzhou 363000, P. R. China\\}
\affiliation{Institute of Meteorological Big Data-Digital Fujian, Minnan Normal University, Zhangzhou, 363000, P. R. China\\}

\begin{abstract}

Numerical computing of the rank of a matrix is a fundamental problem in scientific computation. The datasets generated by the internet often correspond to the analysis of high-dimensional sparse matrices.
Notwithstanding recent advances in the promotion of traditional singular value decomposition (SVD), an efficient estimation algorithm for the rank of a high-dimensional sparse matrix is still lacking.
Inspired by the controllability theory of complex networks, we converted the rank of a matrix into maximum matching computing. Then, we established a fast rank estimation algorithm by using the cavity method, a powerful approximate technique for computing the maximum matching, to estimate the rank of a sparse matrix.
In the merit of the natural low complexity of the cavity method, we showed that the rank of a high-dimensional sparse matrix can be estimated in a much faster way than SVD with high accuracy. Our method offers an efficient pathway to quickly estimate the rank of the high-dimensional sparse matrix when the time cost of computing the rank by SVD is unacceptable.

\end{abstract}

\maketitle

\section{Introduction}
With the development of online social networks, researchers often face complex networks composed of huge numbers of individuals and multiple relationships among them. For the analysis of these large complex networks, we need to convert the network into its corresponding matrix and obtain some characteristics of the original network from its matrix based on traditional matrix theory, such as the page-rank method ~\cite{pagerank2011}, communities detective ~\cite{communities2017}, and some dynamical problem ~\cite{synchronized1998, epidemic2010}.
Rank is one of the most important numerical characteristics of a matrix. At present, a large number of researchers focus on the rank of the special matrix ~\cite{ref5}, low-rank problem ~\cite{ref6,ref7,ref8,ref9}, maximal rank problem ~\cite{ref10}, nullity of graphs ~\cite{ref11, ref12}, and application in robust principal component analysis ~\cite{ref01,ref02,ref03}. The most successful method of rank calculation is the traditional singular value decomposition (SVD), which computes the rank by decomposing the original matrix into singular values and computing the statistical properties of the decomposed matrix. However, the complexity of SVD is the cube of the matrix size (denoted as $N$), which makes the SVD numerically difficult to compute in high-dimensional situations. Therefore, several methods are developed to tackle the complexity problem based on the novel matrix decompositions ~\cite{ref13,ref14}, Monte Carlo simulation ~\cite{ref15}, and multicomputing technologies ~\cite{ref16,ref17,ref18,ref19,ref20}. However, all these methods cannot significantly improve the time complexity.

Benefitting from the development of the control theory of complex networks, we know that the rank of the coupling matrix reflects the exact controllability of sparse complex networks. On the other hand, the structural controllability can be measured by the maximum matching of sparse complex networks. For sparse complex networks, the exact controllability is equivalent to the structural controllability ~\cite{Yuan2013, Liu2011}. The cavity method, a powerful approximation method developed in statistical mechanics ~\cite{zhouhaijun2003, LiuPRL2012}, can be designed to calculate the maximum matching of complex networks. Therefore, the controllability of sparse complex networks builds a bridge between the cavity method and rank computation.

In other words, for an $N$-dimensional sparse matrix, we can convert it into an $N$-node complex network and compute its structural controllability through the cavity method. Due to its sparsity, the structural controllability is equal to the exact controllability, and the rank of the input $N$-dimensional sparse matrix can be approximately estimated. This process, which is a Fast Estimation method for a sparse matrix Rank called FER, can estimate the rank of a high-dimensional sparse matrix much faster than SVD. Therefore, we applied FER to randomly generalized sparse matrices and systematically compared FER with SVD in terms of efficiency, accuracy, and applicability in two typical distributions of nonzero elements of each row (denoted as $\langle k \rangle$). We found that the time cost of FER does not significantly increase during $N$ growing with a constant $\langle k \rangle$, and the results estimated by FER maintained high accuracy, which confirms that FER is an efficient tool for estimating the rank of high-dimensional matrices. We also studied the impact of $\langle k \rangle$ on the time cost and accuracy of FER, and the performance of FER remained very good. Finally, we applied FER to the matrices with the identity of nonzero elements. The efficiency and accuracy of FER were still very high. All the results suggest that FER is a valid access for estimating the rank of a sparse matrix, especially for estimating the rank of a high-dimensional sparse matrix, which is almost unacceptable for computing the rank by SVD while considering the time cost.

\section{Materials and Methods}

FER is based on the development of controllability theory of complex networks. Two existing theoretical frameworks for quantifying the controllability of a complex network are structural controllability theory (SCT) and exact controllability theory (ECT)~\cite{Liu2016}. SCT claims that the structural controllability of any directed network is determined by the maximum matching. The maximum matching can be solved by the cavity method when the network is directed with a structural matrix. The exact controllability obtained by ECT is determined by the maximum multiplicity of eigenvalues of the coupling matrix.
In the sparse situation, ECT is an efficient tool to obtain the controllability of the networks by calculating the rank of the coupling matrix.
When the network is sparse and the weights of links are weakly correlated, the structural controllability and the exact controllability are theoretically equivalent ~\cite{Yuan2013}.
Therefore, computing the rank of a sparse matrix can be converted into a maximum matching problem; then, we can estimate the rank by solving the corresponding coupling equations of the cavity method in an efficient way. This is the core of FER.

\begin{figure}
\begin{center}
\includegraphics[width=\linewidth]{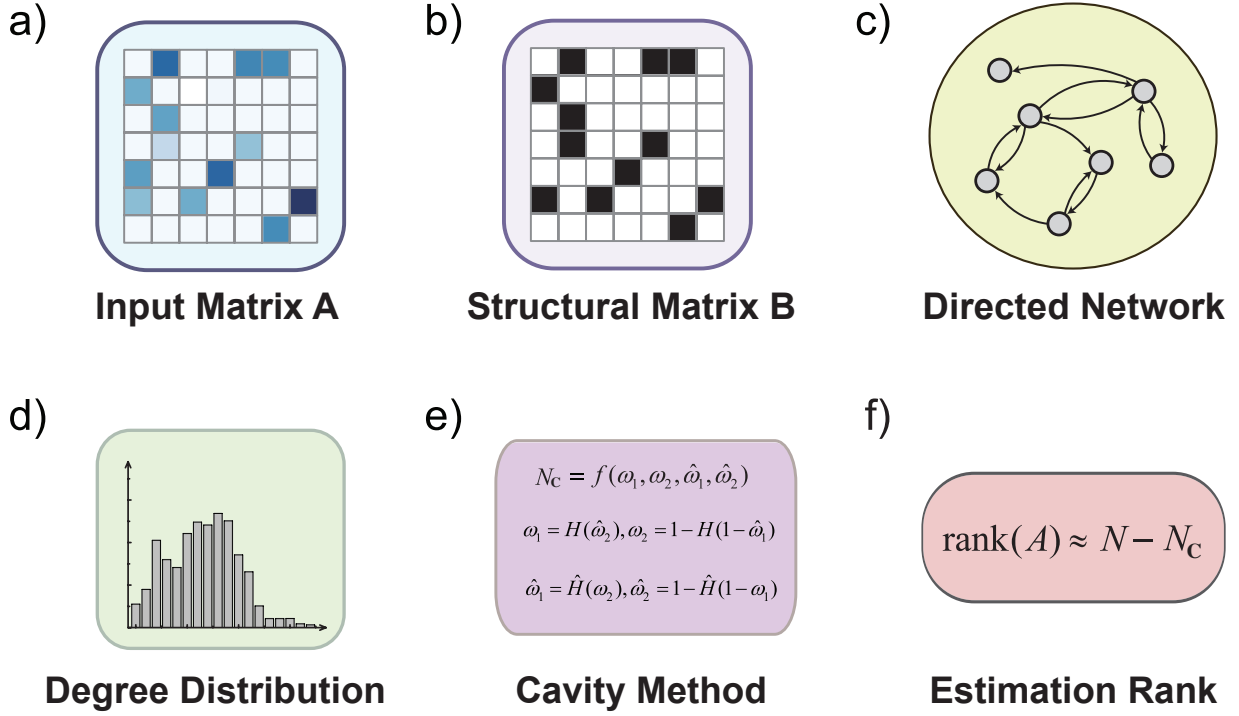}
\caption{{\bf Illustration of the fast estimation algorithm for matrix rank.}
{\bf (a)} Matrix $A$ represents a general sparse matrix as the input matrix, and each grid represents an element in the matrix, in which white grids denote zeros, and darker grids represent the nonzero elements. {\bf (b)} shows the structural matrix of the input matrix $A$. {\bf (c)} transfers the structural matrix $B$ to a directed network. {\bf (d)} shows the statistics for the in-degree and out-degree distributions of the complex network, where the horizontal axis represents the degrees, and the vertical axis represents the relative frequency of the corresponding degree. {\bf (e)} inputs the degree distributions of the network into the coupling equations of the cavity method and solves the values of four coupling parameters. Following eq.~(\ref{eq:nc}), the structural controllability $N_\text{C}$ can be calculated. According to $N_\text{C}$, {\bf (f)} obtains the rank approximation of the input matrix $A$.
}
\label{fig.steps}
\end{center}
\end{figure}

Without loss of generality, we consider an arbitrary sparse input matrix $A$ with weakly correlated nonzero elements, as shown in Fig.~\ref{fig.steps}a, where only the white grids represent the zero elements and the darker color represents the larger value of the nonzero elements. Then, we apply FER to the input matrix $A$, and the procedure of FER can be described as the following five steps:

\textbf{Step 1}. Transfer the input matrix $A$ into a structural matrix $B$, in which the elements can only be $0$ or $1$. $0$s represent the zero elements denoted as white grids, and $1$s represent the nonzero elements denoted as black grids, as shown in Fig.~\ref{fig.steps}a-b;

\textbf{Step 2}. Consider the structural matrix $B$ as a coupling matrix of a complex network and construct a directed network, as shown in Fig.~\ref{fig.steps}b-c;

\textbf{Step 3}. Obtain the in-degree ($P_{in}(k)$) and out-degree ($P_{out}(k)$) distributions of the directed network, where $P_{in(out)}(k)=n_{in(out)}(k)/N$. $n_{in(out)}(k)$ is the number of nodes with the in(out)-degree value $k$ in the whole network, as illustrated in Fig.~\ref{fig.steps}c-d;

\textbf{Step 4}. Calculate the structural controllability ($N_\text{C}$) of the directed network according to the degree distribution by the cavity method~\cite{Liu2011}, which is illustrated in Fig.~\ref{fig.steps}d-e,
\begin{eqnarray}
\label{eq:nc}
N_\text{C} &=& N-\frac{1}{2}[G(\hat{\omega}_2)+G(1-\hat{\omega}_1)-1]+[\hat{G}(\omega_2)+\hat{G}(1-\omega_1)-1]\nonumber \\
&+&\frac{\langle k \rangle}{2}[\hat{\omega}_1(1-\omega_2)+\omega_{1}(1-\hat{\omega}_2)]
\end{eqnarray}
where $G(x)$ and $ \hat{G}(x)$ are ordered by the following equations:
\begin{eqnarray*}
G(x) &=& \sum_{k=0}^{\infty}P_{out}(x)x^{k} \nonumber \text{,}\\
\hat{G}(x) &=& \sum_{k=0}^{\infty}P_{in}(k)x^{k}\text{,}
\end{eqnarray*}
and $\omega_1\text{,}\omega_2\text{,}\hat{\omega}_1\text{,}\hat{\omega}_2$ are the solutions of the following coupling equations:
\begin{eqnarray}
\label{eq:coupling_q}
\omega_1 &=& H(\hat{\omega}_2) \text{,}\nonumber \\
\omega_2 &=& 1-H(1-\hat{\omega}_1) \text{,}\nonumber \\
\hat{\omega}_1 &=& \hat{H}(\omega_2)\text{,} \nonumber \\
\hat{\omega}_2 &=& 1-\hat{H}(1-\omega_1)\text{,}
\end{eqnarray}
in the above equations, the functions of ${H}(*)$ and $\hat{H}(*)$ are shown as:
\begin{eqnarray}
\label{eq:ceq}
H(x) &=& \sum_{k=0}^{\infty}\frac{(k+1)P_{out}(k+1)}{\sum_{k=0}^{\infty}kP_{out}(k)}x^{k}\text{,} \nonumber \\
\hat{H}(x) &=& \sum_{k=0}^{\infty}\frac{(k+1)P_{in}(k+1)}{\sum_{k=0}^{\infty}kP_{in}(k)}x^{k}\text{.}
\end{eqnarray}
According to eq.~(\ref{eq:ceq}), $H(x)$ and $\hat{H}(x)$ can be calculated by the degree distribution ($P_{out}$ and $P_{in}$). In most cases, the degree distribution, a primary statistical property, can be easily obtained from the empirical data as described in step 3.
The above coupling equations are transcendental equations, and the solutions of $\omega_1$, $\omega_2$, $\hat{\omega}_1$, and $\hat{\omega}_2$ can be obtained by numerically solving eq.~(\ref{eq:coupling_q}). Finally, we obtain the structural controllability $N_\text{C}$ from eq.~(\ref{eq:nc}).

\textbf{Step 5}. As the SCT and the ECT are equivalent when the input matrix $A$ is sparse, the structural controllability $N_\text{C}$ is equal to the exact controllability $N-\text{Rank}(A)$. Thus, we can estimate the rank of the input sparse matrix $A$ as illustrated in Fig.~\ref{fig.steps}e-f:

\begin{equation}
\label{eq:rankA}
\text{Rank}(A) \approx  N - N_\text{C}
\end{equation}

 It is worth noting that, first, $N_\text{C}$ can be directly calculated by maximum matching based on SCT~\cite{Liu2016}. The cavity method is an efficient tool based on statistical physics for estimating the maximum matching, which can be obtained just by the degree distribution. That is, the complexity of the FER method is determined by the complexity of the statistics on the degree distribution and the accuracy of the numerical solution. Second, if a matrix contains totally irrelevant element values (every nonzero element is a real random number), $\text{Rank}(A)$ is theoretically equal to $\text{Rank}(B)$ based on the SCT and ECT. However, the assumption is too strict for general cases, which means the result of FER is just an estimation tool for the rank of the input matrix. The correlation strength of nonzero elements in the input matrix indeed affects the accuracy of FER.

\section{Results}
\begin{figure}[htbp]
\begin{center}
\includegraphics[width=\linewidth]{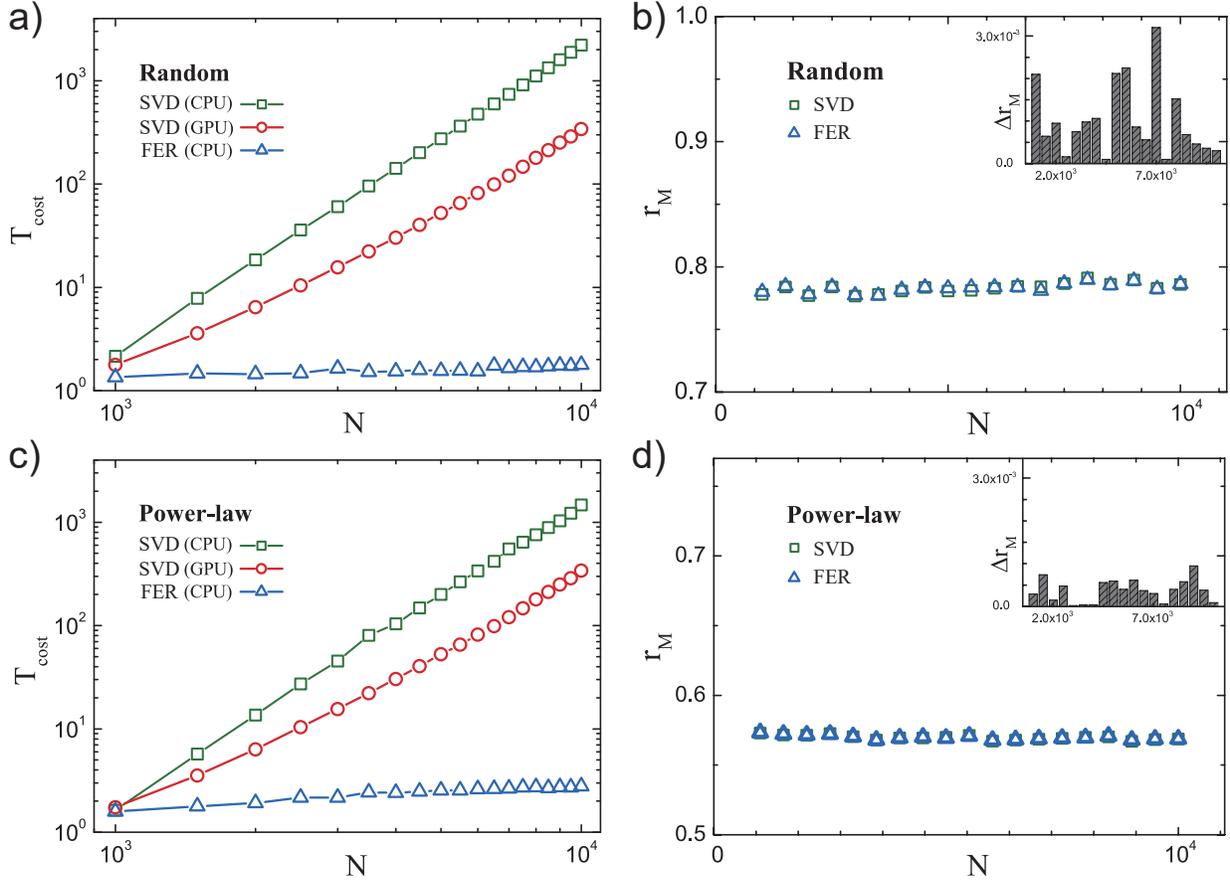}
\caption{
{\bf The impact of $N$ on the efficiency and accuracy of FER.} The changes of $T_{cost}$ to compute the rank of the input matrix using SVD and FER, when input matrix size $N$ increases. $\langle k \rangle$ follows two typical distributions: random distribution {\bf (a)} and power-law distribution ($\gamma=3$) {\bf (c)}. The results of the rank calculated by SVD and the rank estimated by FER in random distribution {\bf (b)} and power-law distribution {\bf (d)} follow eq.~(\ref{r}). The inset figures in {\bf (b)} and {\bf (d)} show the relative errors of FER calculated by eq.~(\ref{rd}). $\langle k \rangle$ is kept at $2$, and all the nonzero elements are random in the generated matrices. The results of $N \leq 5000$ are averaged over 50 independent calculations, and the results of $5000 < N \leq 10000$ are averaged over 20 independent calculations.}
\label{fig2}
\end{center}
\end{figure}

Some comparisons between FER and SVD are exhibited from the efficiency and accuracy aspects in some typical situations.
To analyze the impact of the matrix size ($N$), we generate some matrices randomly with a fixed sparsity, i.e., the average number of nonzero elements in each row ($\langle k \rangle$). The nonzero elements are generated following two typical distributions: random distribution and power-law distribution. Then, we apply FER and SVD to the generated matrices, and the results of comparing the efficiency and the accuracy are shown in Fig.~\ref{fig2}.
 The efficiency of the algorithm is defined by the time cost of solving the task, denoted as $T_{cost}$.
As Fig.~\ref{fig2}a and Fig.~\ref{fig2}c show, if $N$ increases, $T_{cost}$ of SVD increases following its theoretical computational complexity $\mathcal O(N^{3})$. Although we can use a GPU for acceleration, the $T_{cost}$ of SVD increases beyond $\mathcal O(N^{2})$ as $N$ increases.
In contrast, $T_{cost}$ of FER increases very little as $N$ increases, which suggests that its computational complexity is determined by the size of the matrix and the average number of nonzero elements in each row together.
The rank estimated by FER (denoted as $r_\text{M}^{FER}$) and SVD (denoted as $r_\text{M}^{SVD}$) almost overlap, as shown in both Fig.~\ref{fig2}b and Fig.~\ref{fig2}d, which implies that these two methods obtain a similar result no matter how $N$ increases. To explain the high accuracy of FER in more detail, we treat the rank computed by SVD as the ground truth and define the relative error as:
\begin{equation}
\Delta r_\text{M}=\frac{|r_\text{M}^{FER}-r_\text{M}^{SVD}|}{r_\text{M}^{SVD}}\text{,}
\label{rd}
\end{equation}

On the other hand, to compare the FER accuracy as $N$ increases, we define a normalized rank (denoted as $r_{M}$) as the following equation:
\begin{equation}
r_\text{M} \equiv \frac{\text{Rank}(A)}{N}\text{,}
\label{r}
\end{equation}
The inset figures in Fig.~\ref{fig2}b and Fig.~\ref{fig2}d show the relative error $\Delta r_\text{M}$ in random distribution and power-law situations, respectively.
$\Delta r_\text{M}$ are quite small with fluctuations as $N$ increases and remains below $0.003$ and $0.001$ in random  and power-law situations, respectively. The results indicate that FER has good performance in both typical scenarios. When $N$ grows larger, $\Delta r_\text{M}$ has a downward trend in both distributions, which implies that the relative error between FER and SVD should be very small when $N$ is sufficiently large. In summary, for a high-dimensional sparse matrix, we can use FER to obtain an accurate estimation of rank efficiently with a similar accuracy as that obtained by SVD, regardless of the random or power-law distribution.

\begin{figure}[htbp]
\begin{center}
\includegraphics[width=\linewidth]{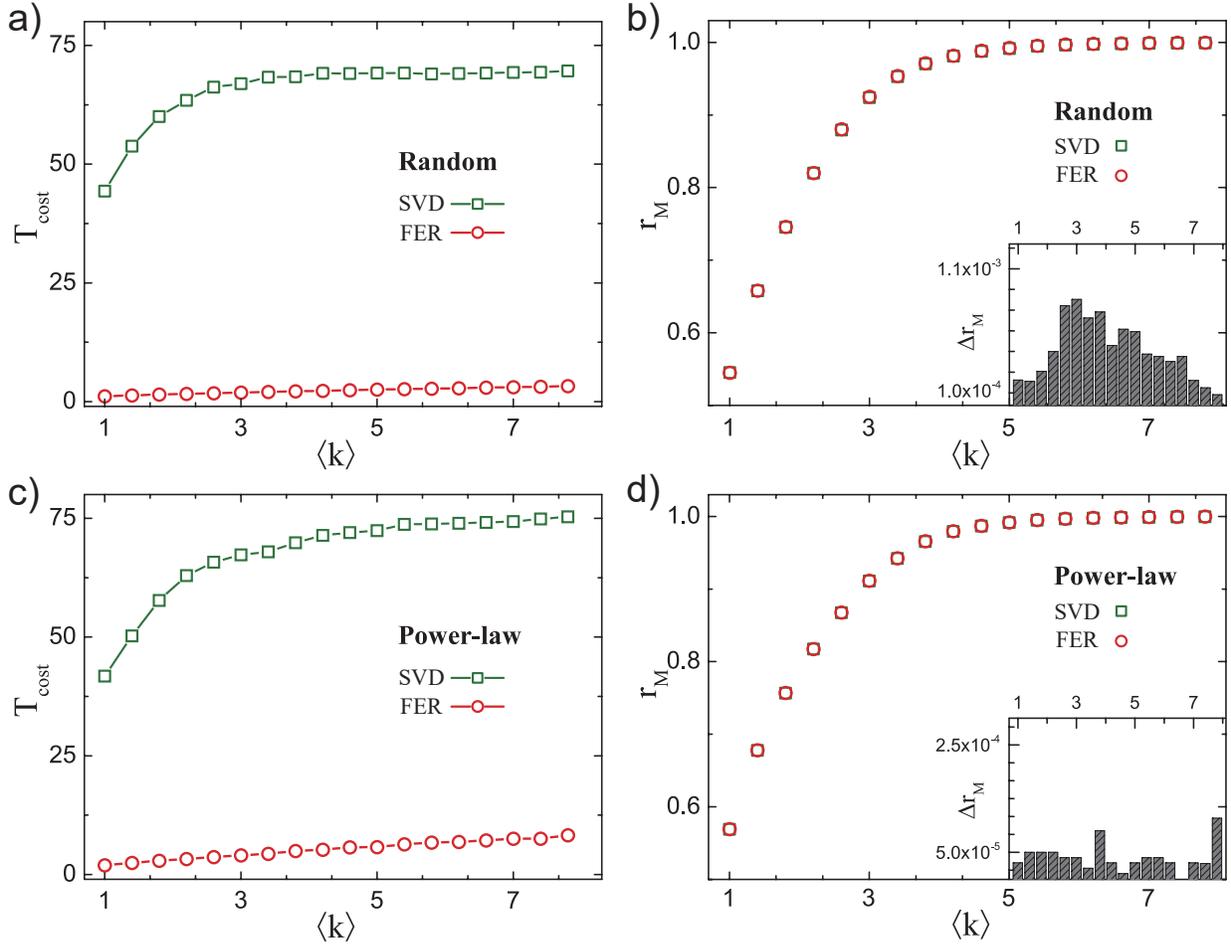}
\caption{
{\bf The impact of $\langle k \rangle$ on FER when $N$ is fixed.} The impact of $\langle k \rangle$ on $T_{cost}^{SVD}$ and $T_{cost}^{FER}$ in random distribution {\bf (a)} and power-law distribution ($\gamma=3$) {\bf (c)}. The comparison between the rank calculated by SVD and FER on generated matrices with random distribution {\bf (b)} or power-law distributions {\bf (d)} for different $\langle k \rangle$. The inset figures in (b,d) show the relative errors of FER versus $\langle k \rangle$. The fixed size of all the simulated networks is $N=3000$, and all the results are averaged over $100$ independent calculations.}
\label{fig3}
\end{center}
\end{figure}

As shown in Fig.~\ref{fig3}, we checked how the sparsity of the input matrix, measured by $\langle k \rangle$, affects the efficiency and accuracy of FER when the matrix size is fixed as $N=3000$.
It is shown that $T_{cost}^{SVD}$ and $T_{cost}^{FER}$ are functions of $\langle k \rangle$ in random situations (Fig.~\ref{fig3}a) and power-law situations (Fig.~\ref{fig3}c). $T_{cost}^{FER}$ is much smaller than $T_{cost}^{SVD}$ in each situation.
In Fig.~\ref{fig3}b and Fig.~\ref{fig3}d, we analyzed the accuracy of FER as $\langle k \rangle$ increased. There are almost no differences between the FER and SVD results, and the scatters in the main figures almost overlap. Then, we consider the relative error of FER, as shown in the inset figures of Fig.~\ref{fig3}b and Fig.~\ref{fig3}d.
If $\langle k \rangle$ increases,

The values of $\Delta r_\text{M}$ are both much smaller in the two situations, which fluctuates obviously in random situations. $\Delta r_\text{M}$ remains almost constant and is smaller than $5 \times 10^{-5}$, in the power-law situation.
In summary, we find that FER is much more efficient than SVD, no matter when $\langle k \rangle$ increases, and the impact of $\langle k \rangle$ on the efficiency and accuracy of FER is quite small in both situations.

\begin{figure}[htbp]
\begin{center}
\includegraphics[width=\linewidth]{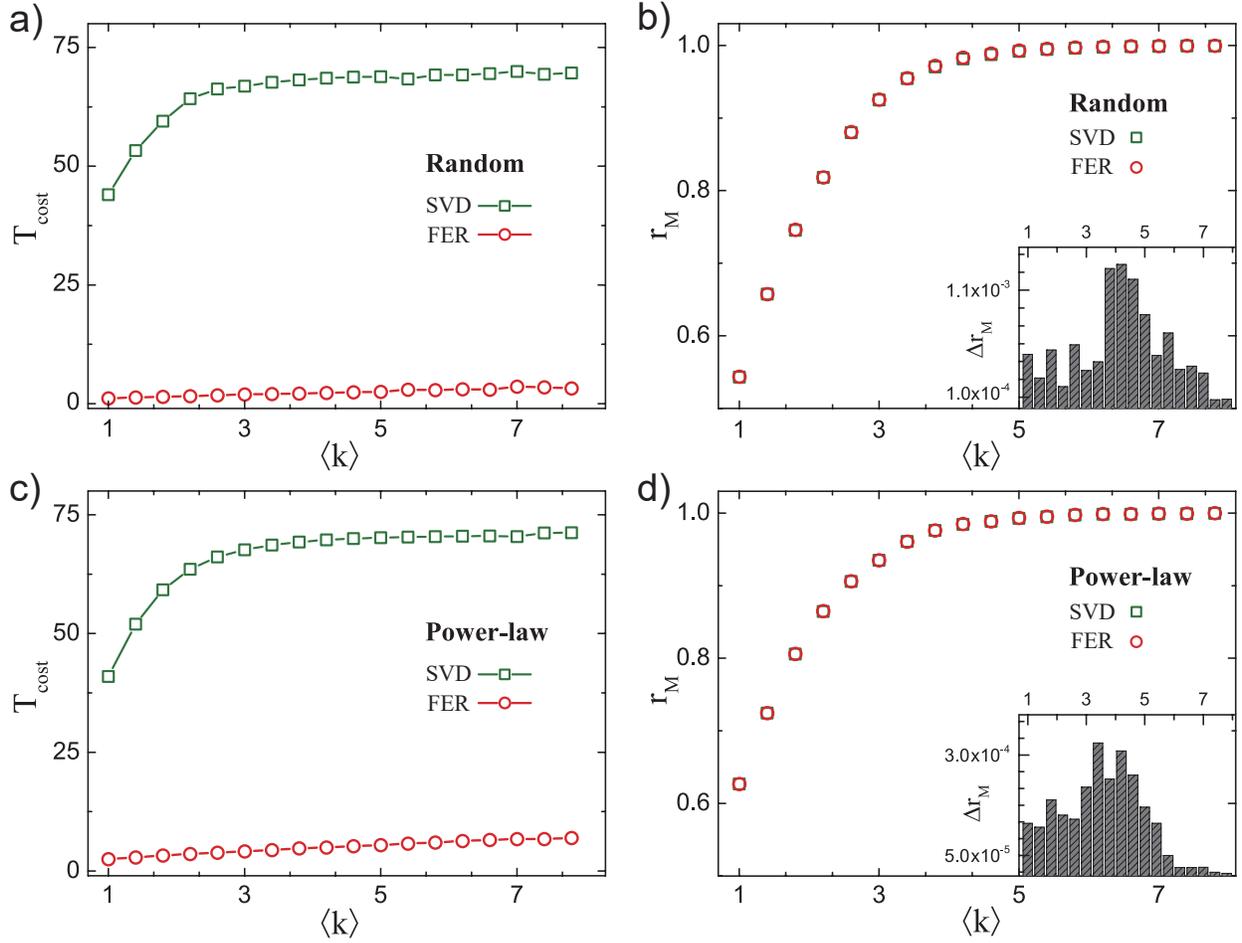}
\caption{
{\bf The efficiency and accuracy of FER when the nonzero elements of the input matrix are all strongly correlated.} $T_{cost}^{SVD}$ and $T_{cost}^{FER}$ on input matrices with random {\bf (a)} or power-law {\bf (c)} distributions for different $\langle k \rangle$. {\bf b,d} the ranks obtained by SVD and FER on generated matrices with random distribution and power-law distributions (d) for different average degrees. All the nonzero elements in the generated matrices are set as $1$. The size of all the simulated networks is $N=3000$, and all the results are averaged over $100$ independent calculations.}
\label{fig4}
\end{center}
\end{figure}

FER works only if all the nonzero elements in the sparse matrix are uncorrelated. However, there are many relevant elements in the real data, which means that errors are unavoidable if the nonzero elements are correlated. Thus, we discuss whether the result obtained approximately by FER is acceptable when the nonzero elements of the input matrix are correlated. In Fig.~\ref{fig4}, we consider an extreme case where all the nonzero elements in the input matrix are identical (set as $1$), which means that all the nonzero elements are strongly correlated.
The strong correlation has a negative effect on $T_{cost}$ in both random situations (Fig.~\ref{fig4}a) and power-law situations (Fig.~\ref{fig4}c). $T_{cost}^{FER}$ is still much smaller than $T_{cost}^{SVD}$. However, the results shown in Fig.~\ref{fig4}b and Fig.~\ref{fig4}d indicate that the accuracy of FER has a significant decline compared with Fig.~\ref{fig3}b and Fig.~\ref{fig3}d. This means that the correlation of the input matrix does affect the accuracy of FER, which agrees with the limitation of structural controllability as well as the cavity method. Although the accuracy of FER has decreased, we can also learn from the inset figures that the relative error of FER is still very small in both situations. Especially when $\langle k \rangle$ increases over $4$, $\Delta r_\text{M}$ has an obvious descent.
In other words, even though the nonzero elements in the sparse matrix are strongly correlated, the performance of FER is also acceptable in terms of efficiency and accuracy. The robustness of FER suggests its effectiveness in estimating the rank of a more general matrix extracted from the empirical data set.

\section{Discussion}
In summary, we utilized the cavity method to estimate the maximum matching. Based on controllability theory of complex networks, we know that the rank of the matrix is theoretically equal to the maximum matching of the network when the network is sparse and the weights of links are weakly correlated. Then, we established an efficient estimation tool for analyzing the rank of a high-dimensional sparse matrix by the cavity method, which is called FER. We discussed the impact of the input matrix size ($N$), the sparsity of the matrix (measured by $\langle k \rangle$), and the correlation of nonzero elements on the efficiency (measured by $T_{cost}$) and accuracy (measured by $\Delta r_\text{M}$) of FER in random situations and power-law situations.
We found that FER has remarkable performance in terms of both efficiency and accuracy in random distribution and power-law distribution. Although the characteristics of nonzero elements affect the results, FER can still be applied to most sparse matrices to estimate their rank with fast speed and high accuracy. It can significantly outperform SVD in terms of the time cost and has a similar accuracy to SVD. Therefore, FER provides an efficient and accurate method for estimating the rank of a sparse matrix. Especially for dealing with a large real network by some algorithms with its matrix rank, FER can do a good job to estimate the rank directly by its degree distribution obtained from the raw data, while SVD is inapplicable due to its excessive time cost. Furthermore, in some special situations, where only the structural information of a social network can be detected, such as degree distribution or partially missing degree distribution, FER is still applicable to estimate the rank of its corresponding matrix. This means that FER can potentially be used in some algorithms designed for incomplete data or data polluted by interference noise.

\section*{Acknowledgements}
We thank Professor Wen-Xu Wang for valuable suggestions. This work is supported by the National Natural Science Foundation of China (Grant Nos. 61703136 and 61672206), the Natural Science Foundation of Hebei (Grant Nos. F2020205012 and F2017205064), and the Youth Excellent Talents Project of Hebei Education Department (Grant No. BJ2020035).

\end{document}